\documentclass[12pt,twoside]{article}

\usepackage{times,amsmath,amsthm,amssymb,bm}
\usepackage{type1cm}

\numberwithin{equation}{section}

\theoremstyle{plain}
\newtheorem{thm}{Theorem}[section]
\newtheorem{lem}[thm]{Lemma}
\newtheorem{prop}[thm]{Proposition}
\newtheorem{defn}[thm]{Definition}

\newtheorem*{prob1}{Linearized Inverse Problem of the real Drift (LIPD)}

\newtheorem{rem}{Remark}[section]

\newcommand{\memo}[1]{}

\begin{document}

\title{\vspace*{-5.5ex}
\bf{Application of Microlocal Analysis to an Inverse Problem Arising from Financial Markets }\vspace*{-2.5ex}}
\author{{\normalsize by}\\[0.5ex]  Shin-ichi Doi and Yasushi Ota}
\date{\normalsize September 09, 2014}
\maketitle

\footnote{ 
{\bf 2000 Mathematics Subject Classifications:} 
Primary 
35R30; 
Secondary 
35K08. 
}

\footnote{ 
{\bf Keywords and Phrases:} 
Inverse problem;
Microlocal Analysis;
FBI transform.
}

\vspace*{-11ex}

\begin{abstract}\vspace*{0ex}
\begin{center}\hspace*{-1.5ex}\begin{minipage}{62ex}\hspace{2.5ex}
One of the most interesting problems discerned when applying the Black--Scholes model to financial derivatives, is reconciling the deviation between expected and observed values.
In our recent work, we derived a new model based on the Black--Scholes model and formulated a new mathematical approach to an inverse problem in financial markets. 
In this paper, we apply microlocal analysis to prove a uniqueness of the solution to our inverse problem. 
While microlocal analysis is used for various models in physics and engineering, this is the first attempt to apply it to a model in financial markets.
First, we explain our model, which is a type of arbitrage model.
Next we illustrate our new mathematical approach, and then for space-dependent real drift, we obtain stable linearization and an integral equation. Finally, by applying microlocal analysis to the integral equation, we prove our uniqueness of the solution to our new mathematical model in financial markets.
\end{minipage}
\end{center}\vspace*{6ex}
\end{abstract}

\section{Introduction}

Financial derivatives are contracts wherein payment is derived from an underlying asset such as a stock, bond, commodity, interest, or exchange rate. 
An underlying asset $S_t$ at time $t$ is modeled by the following stochastic differential equation:
\begin{equation*}
dS_t =  \mu(t,S_t) S_t dt + \sigma(t,S_t)S_t dW_t,
\label{sde}
\end{equation*}
where the process $W_t$ is Brownian motion. The parameters $\mu(t,S)$ and $\sigma(t,S)$ are called the real drift and the local volatility of the underlying asset, respectively. 

Black and Scholes\cite{B-S} first discovered how to construct a dynamic portfolio $\Pi_t$ of a derivative security and the underlying asset.
Their approach is developed in probability theory, and the hedging, and pricing theory of the derivative security is established as mathematical finance.
By Ito$'$s lemma, the stochastic behavior of the derivative security $u(t, S)$ is governed by the following stochastic differential equation:
\begin{equation*}
du = \biggl(\frac{\partial u}{\partial t}+\mu(t, S)S \frac{\partial u}{\partial S} + \frac{1}{2} \sigma(t, S)^2\frac{\partial^2 u}{\partial S   ^2}\biggr) dt + \sigma(t, S)S \frac{\partial u}{\partial S} dW.
\label{spde}
\end{equation*}
In the absence of arbitrage opportunities, the instantaneous return of this portfolio must be equal to the interest rate $r>0$, i.e., the return on a riskless asset such as a bank deposit.
Therefore, this equality takes the form of the following partial differential equation:
\begin{equation}
\frac{\partial u}{\partial t}+\frac{1}{2}{{\sigma}(t,S)}^{2}S^{2}\frac{\partial^{2}u}{\partial S^{2}}+ (r-\delta) S \frac{\partial u}{\partial S}-ru = 0,    
\label{bseq}
\end{equation}
where $r$ and the divided rate $\delta$ are the known constants.

Their approach provides a useful, simple method of pricing inclusive of financial derivatives, risk premium, and default probability estimation under the assumption that the risky asset is log-normally distributed.
However, the theoretical prices of options with different strike prices as  calculated by the Black--Scholes model differ from real market prices. 
Specifically, when we apply the Black--Scholes model to default probability estimation, we must be careful of the deviation that arises between expected and observed values.
Merton\cite{merton} has formulated a default probability estimation using a model based on \cite{B-S} by considering the value of the firm instead of its stock, the firm$'$s debt instead of strike price, and its equity instead of option price and Boness\cite{Boness} has derived the formulation of it by another method.    
However, as shown in deriving the Black-Scholes model (see \cite{B-S}), under the no-arbitrage property of the financial market, the real drift $\mu$ does not enter equation (\ref{bseq}).
In \cite{myy}, taking this into account, we have derived the following new model by using $A_t$ instead of $S_t$:
\begin{equation}
\frac{\partial u}{\partial t}+\frac{1}{2}{{\sigma}(t,A)}^{2}A^{2}\frac{\partial^{2}u}{\partial A^{2}}+ \mu(t,A) A \frac{\partial u}{\partial A}-ru = 0.    
\label{maineq}
\end{equation}
Moreover, in \cite{myy} we have established an inverse problem to reconstruct the real drift from the observable data, but only an binary option case. 
In Korolev, Kubo and Yagola\cite{kubo}, they reconstructed the unknown drift in our new model.

In this paper, we prove the uniqueness of the solution to an inverse problem with respect to the real drift by applying microlocal analysis. 
To give a brief description of our problem, we build upon the method in \cite{B-I-V}. 
In \cite{B-I-V}, they used the standard linearization method with an option pricing inverse problem and derived the partial differential equation with the constant coefficient $\sigma_0^2, \delta, r$ for the linear principal part $V$.
Since a change of variables means 
this equation is reduced to the heat equation with the right-hand side 
$w(\tau,y) f(y) $, they wrote the well--known integral representation for the solution $W$ to that heat equation with a suitable initial condition as follows:
\begin{equation}
\displaystyle  W(\tau, x) = \displaystyle \int_{\textbf{R}} \int_0^{\tau} \frac{1}{\sqrt{2 \pi (\tau - \theta )\sigma_0^2}} e^{-\frac{|x-y|^2}{2 \sigma_0^2 (\tau - \theta) }} w(\theta, y) f(y) d\theta dy, 
\label{BIVeq}
\end{equation}
where $f$ is a small perturbation of constant $\sigma_0$, $w(\tau, y)$ is represented by 
\begin{equation*}
w(\tau, y) = \frac{s^*}{\sqrt{2\pi \tau \sigma_0^2}} e^{- \frac{|y|^2}{ 2 \sigma_0^2 \tau}},
\label{bivw}
\end{equation*}
Here $\tau = T - t$, $y=\log K/s^*$, $K$ is a strike price at the maturity date $T$ and $s^*$ is market price of the stock at a current time $t^*$.

For the above equation, they applied the Laplace transform to exactly evaluate an integral with respect to time. As a result, they derived the integral equation for $f$ that takes the following form
\begin{equation}
\displaystyle V(\tau, x) = \int_{\textbf{R}} B(x, y; \tau) f(y) dy   
\label{iebiv}
\end{equation}
with the kernel 
\begin{equation*}
B(x,y; \tau) = \displaystyle \frac{s^*}{\sigma_0^2 \sqrt{\pi}} \int_{\frac{|x-y|+|y|}{\sigma_0\sqrt{2 \tau}}}^{\infty} e^{-\theta^2} d\theta 
\end{equation*}
given by the error function, and thus proved the uniqueness for the linearized inverse problem.
In our case, since our principal linear part $W$ which is derived in the same manner as \cite{B-I-V} has the following form
\begin{equation}
W(\tau, x) = \displaystyle \int_{\textbf{R}} \int_{0}^{\tau}  \frac{1}{\sqrt{4 \pi (\tau - s) \sigma_0^2}} e^{-\frac{|y-x|^2}{4 (\tau - \theta) \sigma_0^2 }} w(\theta,y) f(y) d\theta dy, 
\label{ourieq}
\end{equation}
where $w(\theta, y)$ takes the following form
\begin{equation*}
w(\tau, y) = \displaystyle \int_0^{\infty} \frac{1}{\sqrt{4 \pi \tau \sigma_0^2 }} e^{-\frac{|x-y|^2}{4 \tau \sigma_0^2} } dx.
\end{equation*}
Therefore we are unable to derive an integral equation by the Laplace transform as in (\ref{iebiv}); that is, in our case $w(\tau, y)$ is not a Gauss function but an error function.
In this paper, taking this into account, we shall prove the uniqueness of the solution to the inverse problem of the real trend by applying the Fourier-Bros-Iagolnitzer (for short, FBI) transform to (\ref{ourieq}).

The paper is divided into six sections. 
In Section 2, we illustrate the Linearized Inverse Problem of the real Drift (LIPD). 
The main theorem is stated in Section 3 wherein we provide the outline of main theorem. 
In Section 4, we summarize basic facts concerning the FBI transform which plays an essential role in the proof of our main theorem.
The proof of main theorem is proved in Section 5 and some mathematical results used in Section 5 are proved in Section 6.

\section{Inverse problem of the real drift}

In \cite{myy}, we have derived a new arbitrage model and formulated an inverse option pricing problem for a reconstruction of a real drift in the binary option case. In this section, we explain how to formulate an inverse problem of our new arbitrage model and reconstruct the real drift.

Here, we consider the following problem wherein the local volatility $\sigma(t, A)$ is a positive constant $\sigma_0 > 0$ and the real drift $\mu(t, A)$ is a time-independent in our new equation (\ref{maineq}) with a suitable condition: 
\begin{equation}
u(t, A)|_{t=T} = \max \{ A-D, 0 \}
\label{inicondi}
\end{equation}
where $D$ is a price of the firm$'$s debt at the maturity date $T$.

By the following changes of variables and substitutions 
\begin{align}
y = \log \displaystyle \frac{A}{D},& \hspace{5ex} \tau = T - t, \notag \\
\mbox{} \\
\mu(y) = \mu(De^y), & \hspace{5ex} U(\tau, y) = u(T-\tau, De^y)/D,  \notag 
\end{align}
the equation (\ref{maineq}) and the initial data can be transformed into the following form:
\begin{equation}
\left\{
\begin{array}{cl}
\displaystyle\frac{\partial U}{\partial \tau}  = \frac{1}{2} \sigma_0^2 \frac{\partial^2 U}{\partial y^2} - \Bigl(\displaystyle\frac{1}{2} \sigma_0^2 - \mu(y) \Bigr)\frac{\partial U}{\partial y} - r U \hspace{3ex} (y, \tau) \in \textbf{R} \times(0, \tau^*), \\
\mbox{} \\
U(\tau, y)|_{\tau=0}  = \max \{ e^y - 1, 0 \} \hspace{10ex}  y \in \textbf{R}, 
\end{array}
\right.
\label{mueq1}
\end{equation}
\begin{eqnarray}
U(\tau^*, y) = U^*(y)   \hspace{20ex}  y \in \omega \subseteq \textbf{R}, 
\label{mueq1add}
\end{eqnarray}

\noindent where $\tau^* = T - t^*>0 $, $t^*$ is the current time and $\omega$ is an interval of $\textbf{R}$.

Here we define that the inverse problem of the real drift (\ref{mueq1}) and (\ref{mueq1add}) seeks $\mu(y)$ from the given $U^*(y)$. However, since this inverse problem is nonlinear, difficulties arise with the uniqueness and existence of the solution.
Therefore, we will formulate the inverse problem of the real drift by means of the linearization method in \cite{B-I} and \cite{B-I-V}.

To linearize around the constant coefficient $\mu_0$, we assume that
\begin{eqnarray*}
\mu(y)=\mu_0+f(y),
\end{eqnarray*}
where $f(y)$ denotes a small perturbation. Thus, we observe
\begin{eqnarray*}
U=U_0+V+\nu,
\end{eqnarray*}
where $U_0$ solves the Cauchy problem (\ref{mueq1}) with $\mu(y)\equiv \mu_0$, $\nu $ is quadratically small with respect to $f$, and $V$ is the principal part of the perturbed solution $U$.
Substituting this into the expression for $u$ and neglecting terms of higher order with respect to $f$, we reach the linearized inverse problem of the real drift.

\vspace{2ex}
\begin{prob1}{}
The parameters $\tau^*$, $\mu_0$, $\sigma_0$, and $r$ are given. From the option price $V^*(y)=\left\{U^*(y)-U_0(\tau^*,y)\right\}$, identify the perturbation $f(y)$ satisfying
\begin{eqnarray}
\left\{
  \begin{array}{c}
    \displaystyle\frac{\partial V}{\partial \tau}-\frac{1}{2}\sigma_0^2\frac{\partial^2V}{\partial y^2}+\left(\frac{1}{2}\sigma_0^2-\mu_0 \right)\frac{\partial V}{\partial y}+rV=\frac{\partial u_0}{\partial y}f(y),\vspace{5mm}   \\
V(\tau,y)|_{\tau=0}=0,   \\
  \end{array}
\right.
\label{mueq2}
\end{eqnarray}
\begin{eqnarray}
V(\tau^*,y)=V^*(y).
\label{mueq2add}
\end{eqnarray}
\end{prob1}

\section{Main results}

In this section we prove the uniqueness of the solution to LIPD by using microlocal analysis.
Before describing the main theorem, we shall transform equation (\ref{mueq2}) into simple form and derive a Fredholm--type integral equation.  

We set 
\begin{eqnarray*}
a_0 = \displaystyle \frac{\sigma_0^2 - 2 \mu_0}{2 \sigma_0^2}, \hspace{3ex} b_0 = r + \frac{1}{2} \sigma_0^2 a_0^2, \\
H_a = - \left(\frac{\partial}{\partial y} - a \right)^2 \hspace{3ex} (a= a_0-1), 
\end{eqnarray*}
then (\ref{mueq2}) can be rewritten as 
\begin{eqnarray}
\left\{
\begin{array}{cl}
\displaystyle \left( \frac{\partial}{\partial \tau} + \frac{1}{2}\sigma_0^2 H_a \right) v(\tau, y) = f(y)w(\tau, y) \hspace{3ex} (\tau, y) \in (0, \tau^*) \times\textbf{R}, \\
\mbox{} \\
v(\tau, y)|_{\tau=0}  = 0 \hspace{25ex}  y \in \textbf{R}, 
\end{array}
\right.
\label{mueq3}
\end{eqnarray}
where, $v(\tau, y) = e^{-y + b_0 \tau} V(\tau, y)$ and $w(\tau, y)$ is the following form
\begin{equation}
w(\tau, y) = e^{-y + b_{0}\tau} \frac{\partial U_0}{\partial y}.
\end{equation}
Here $w$ is the solution of the following problem
\begin{eqnarray}
\left\{
\begin{array}{cl}
\displaystyle \left( \frac{\partial}{\partial \tau} + \frac{1}{2}\sigma_0^2 H_a \right) w(\tau, y) = 0 \hspace{10ex} (\tau, y) \in (0, \tau^*) \times\textbf{R}, \\
\mbox{} \\
w(\tau, y)|_{\tau=0}  = 1_{[0, \infty]}(x) \hspace{25ex}  y \in \textbf{R}, 
\end{array}
\right.
\label{mueq-w}
\end{eqnarray}
Now by setting
\begin{equation*}
\tilde v(\tau, y) = v\left(\frac{2}{\sigma_0} \tau, y\right) , \hspace{1ex} \tilde w(\tau, y) = w\left(\frac{2}{\sigma_0} \tau, y\right) 
\end{equation*}
and 
\begin{equation*}
\tilde{f}(y) = \frac{2}{\sigma_0}f(y), \hspace{1ex} \tilde{\tau^*} = \frac{\sigma_0}{2}\tau^*,
\end{equation*}
(\ref{mueq3}) and (\ref{mueq-w}) can be rewritten the following simple form

\begin{eqnarray}
\left\{
\begin{array}{cl}
\displaystyle \left( \frac{\partial}{\partial \tau} +  H_a \right) \tilde{v}(\tau, y) = \tilde{f}(y)\tilde{w}(\tau, y) \hspace{3ex} (\tau, y) \in (0, \tilde{\tau}^*) \times\textbf{R}, \\
\mbox{} \\
\tilde{v}(\tau, y)|_{\tau=0}  = 0 \hspace{25ex}  y \in \textbf{R}, 
\end{array}
\right.
\label{mueq3-new}
\end{eqnarray}
and
\begin{eqnarray}
\left\{
\begin{array}{cl}
\displaystyle \left( \frac{\partial}{\partial \tau} +  H_a \right) \tilde{w}(\tau, y) = 0 \hspace{10ex} (\tau, y) \in (0, \tilde{\tau}^*) \times\textbf{R}, \\
\mbox{} \\
\tilde{w}(\tau, y)|_{\tau=0}  = 1_{[0, \infty]}(x) \hspace{25ex}  y \in \textbf{R}.
\end{array}
\right.
\label{mueq-w-new}
\end{eqnarray}
From now we consider problems (\ref{mueq3-new}) and (\ref{mueq-w-new}) wherein $\tilde{v}$, $\tilde{w}$, $\tilde{f}$ and $\tilde{\tau}^*$ are rewritten as $v$, $w$, $f$ and $\tau$ respectively, if there is no confusion. 

By the well--known representation of the solution to the Cauchy problem (\ref{mueq3-new}), we have the following Fredholm--type integral equation: 
\begin{eqnarray}
v(\tau^*, x) = \displaystyle\int_{0}^{\tau^*} U_a(\tau^* - s) [w(s, \cdot) f(\cdot)](y) ds.
\label{inteq1}
\end{eqnarray}
Here 
\begin{eqnarray*}
(U_a(\tau) \varphi)(y) = \displaystyle \int_{\mathbf{R}} K_a (\tau, y-x) \varphi(x) dx, 
\end{eqnarray*}
where 
\begin{eqnarray*}
K_a(\tau, y) = \displaystyle \frac{1}{\sqrt{4 \pi \tau}} e^{-\frac{|y|^2}{4 \tau}+ ay}
\end{eqnarray*}
and $w(\tau, x)$ is represented by the following form:
\begin{align}
w(\tau, x) := &(U_a(\tau) H_+)(x) \notag \\                                                                                                                                                                                                               = & \displaystyle \int_0^{\infty} \frac{1}{\sqrt{4 \pi \tau  }} e^{-\frac{|x-y|^2}{4 \tau } + a(x-y)} dx \notag \\ 
=& \displaystyle \frac{1}{\sqrt{\pi}} e^{\tau a^2} \int_{-\infty}^{\frac{x-2\tau a}{\sqrt{4 \tau}}} e^{-\theta^2} d\theta,
\label{w}
\end{align}
where $H_+(x) = 1_{[0, \infty]}(x)$.

We will describe the results for LIPD in the following theorem.
\begin{thm}\label{mainthm}
Let $\tau^* > 0$ and $f(y) \in L^2(\textbf{R}).$ Assume that $ supp f \subset [-L, \infty)$ with some $L \ge 0$. 
Then a solution f(y) to the integral equation (\ref{inteq1}) and hence to the inverse problem of the real drift (\ref{mueq2}) and (\ref{mueq2add}) is unique.  
\end{thm}
\noindent \textit{Outline of Proof.}
To prove the claim of Theorem \ref{mainthm}, it suffices to prove $f = 0$ under the assumption that the left--hand side of (\ref{inteq1}) is zero.

Now we assume that $v(\tau^*, y)$ is zero, that is, 
\begin{equation}
\displaystyle\int_{0}^{\tau^*} U_a(\tau^* - s) [w(s, \cdot) f(\cdot)](y) ds = 0.
\label{Ieq0}
\end{equation}
To prove that $f$ is zero, we will show that there exist $\delta > 0$ such that 
\begin{equation}
||T f ||_{L^2([-L_0, \infty) \times \{|\xi| \ge 2 \})} = O(e^{-\frac{\delta}{h}}) 
\label{outline}
\end{equation}
where $L_0= L+1$.
Here, $Tf$ is called the \textit{Fourier-Bros-Iagolnitzer} (for short, FBI) transform of $f$. 
By Definition \ref{analytic} in the next section, since the estimation of (\ref{outline}) lead us to the following assertion 
\begin{equation*}
[-L_0, \infty) \times \{ |\xi| \ge 2 \} \cap {\rm WF}_a(f) = \emptyset,
\end{equation*}
where ${\rm WF}_a$ is called the analytic wave front set of $f$, we obtain that $f$ is real analytic in $(-L_0, \infty)$.
Moreover, since $f=0$ in $(-L_0, -L)$ by the assumption, we will be able to conclude that $f$ is identically zero on $\textbf{R}$.

Therefore, to prove the claim of Theorem \ref{mainthm}, it suffices to derive the estimation (\ref{outline}) under the assumption that the left--hand side of (\ref{inteq1}) is zero.

\section{Properties of FBI transform}

In this section, we summarize basic facts concerning the FBI transform (see \cite{Mar}).

\begin{defn}
For $u\in \mathcal{S}'(\mathbf{R}^n)$, the FBI transform of $u$, $Tu$, is defined as 
\begin{align*}
Tu(x, \xi; h)= (2\pi h)^{-\frac{n}{2}}(\pi h)^{-\frac{n}{4}} \displaystyle \int_{\mathbf{R}^n} e^{i(x-y)\cdot\xi/h-(x-y)^2/2h} u(y)dy,\ (x,\xi)\in \mathbf{R}^{2n}. 
\end{align*}
Here the integral is in the sense of distributions, and $h>0$ is a parameter.
(Parameter $h$ is often omitted if there is no confusion.) 
\label{defofFBI}
\end{defn}

\begin{rem}
For $u\in L^2(\mathbf{R}^n)$, $\|Tu(\ \cdot \ , \ \cdot\ \ ;h)\|_{L^2(\mathbf{R}^{2n})}=\|u\|_{L^2(\mathbf{R}^n)}$.
\end{rem}


\begin{defn}\label{analytic}
A distribution $u\in \mathcal{S}'(\mathbf{R}^n)$
is called analytic at $(x_0,\xi_0)\in \mathbf{R}^n\times(\mathbf{R}^n\setminus\{0\})$
if there exist $\delta>0$ and a neighborhood $V$ of $(x_0,\xi_0)$ such that 
\begin{equation}
\|Tu(\ \cdot \ , \ \cdot\ \ ;h)\|_{L^2(V)}=O(e^{-\delta/h}) \text{ as $h\to +0$}.
\label{ex-small}
\end{equation}
The analytic wave front set of $u$, $WF_a(u)$, is the set of 
all $(x_0,\xi_0)\in \mathbf{R}^n\times(\mathbf{R}^n\setminus\{0\})$
at which $u$ is not analytic.
\end{defn}

\begin{rem}\label{a-wave-front}
(i) We can replace $\|Tu\|_{L^2(V)}$ by $\|Tu\|_{L^{\infty}(V)}$ in the definition above.

(ii) If $(x,\xi)\in WF_a(u)$, then $(x,t\xi)\in WF_a(u)$ for every $t>0$.

(iii) $u$ is real-analytic near $x_0$ if and only if $\left(\{x_0\}\times(\mathbf{R}^n\setminus\{0\})\right)\cap WF_a(u)=\emptyset$.
\end{rem}

\begin{defn}
Put $\langle \xi \rangle=\sqrt{1+|\xi|^2}\ (\xi\in \mathbf{R}^n)$.
Let $m\in \mathbf{R}$.
The symbol class $S_{2n}(\langle \xi \rangle ^m)$ is the set of all $p=p(x,\xi;h):\mathbf{R}_{x,\xi}^{2n}\times (0,1)\to\mathbf{C}$
satisfying the following conditions:

(i) For each $h$, $p(\cdot\ ;h)\in C^\infty(\mathbf{R}^{2n})$.

(ii)  For every $\alpha,\beta\in \mathbf{N}_0^n$, there exists $C>0$ such that
$$
|\partial_x^\alpha\partial_\xi^\beta p(x,\xi;h)|\le C \langle \xi \rangle ^m \quad 
\text{for every $(x,\xi;h)\in \mathbf{R}^{2n}\times (0,1)$}.
$$
\end{defn}

\begin{defn}
For $u\in \mathcal{S}(\mathbf{R}^n)$, $F_h u $ is defined as follows
\begin{equation}
F_h u(\xi) = \frac{1}{(2\pi h)^{\frac{d}{2}}} \int_{\mathbf{R}^d} e^{-ix\xi/h} u(x) dx,
\end{equation}
where $\xi \in \mathbf{R}^d$ and $x \xi$ stand for the scalar product of $x$ and $\xi$.
\end{defn}

\begin{defn}
Let $t\in [0,1]$ be fixed.
For $p\in S_{2n}(\langle \xi \rangle ^m)$ and $u\in \mathcal{S}(\mathbf{R}^n)$, 
${\rm Op}_h^t(p)u \in \mathcal{S}(\mathbf{R}^n)$ is defined as follows
\begin{equation*}
{\rm Op}_h^t(p)u(x)=
\frac{1}{(2\pi h)^n}\displaystyle \int_{\mathbf{R}^n}\int_{\mathbf{R}^n} e^{i(x-y)\cdot\xi/h} p((1-t)x+ty,\xi;h)u(y)dyd\xi.
\end{equation*}
The operator ${\rm Op}_h^t(p)$ can be extended as a continuous operator in $\mathcal{S}'(\mathbf{R}^n)$.
\end{defn}

\begin{rem}
In this paper, we use only the case $t=1$.
\end{rem}

\begin{defn}
Let $a>0$, and set $\Sigma_a=\{x\in \mathbf{C}^n;\ |\mathrm{Im}\, x|<a\}$.
The symbol class $S_{2n}^{hol}(\langle \xi \rangle ^m,\Sigma_a)$ is the set of all $p=p(x,\xi;h):\Sigma_a\times \mathbf{R}^{n}\times (0,1)\to\mathbf{C}$
satisfying the following conditions:

(i) $p=p(x,\xi;h)$ is $C^\infty$ in $(x,\xi)\in \Sigma_a\times\mathbf{R}^n$ and holomorphic in $x\in \Sigma_a$.

(ii)  For every $\alpha,\beta\in \mathbf{N}_0^n$, there exists $C>0$ such that
$$
|\partial_x^\alpha\partial_\xi^\beta p(x,\xi;h)|\le C\langle \xi \rangle ^m\quad \text{for every $(x,\xi;h)\in \Sigma_a\times\mathbf{R}^n\times (0,1)$}.
$$
\end{defn}

\begin{prop}\label{prop-f-supp}
Let $F_1, F_2$ be a closed sets in $ \mathbf{R}^d $ satisfying 
\begin{equation}
{\rm dist}(F_1, F_2) = \sigma_{1} > 0.
\label{prop5.2}
\end{equation}
Then for $u \in L^2(\mathbf{R}^d)$ with ${\rm supp} \ u \subset F_1$, 
there exists $\delta > 0$ such that 
\begin{equation}
||T u (\ \cdot \ , \ \cdot \ ; h)||^2_{L^2(F_2 \times \mathbf{R}^d )} \le 2^{\frac{d}{2}}e^{- \frac{\delta}{h}} ||u||_{L^2(\mathbf{R}^d)}.
\end{equation}
\begin{proof}
By definition of the FBI transform, we have 
\begin{align}
||T u (\ \cdot \ , \ \cdot \ ; h)||^2_{L^2(F_2 \times \mathbf{R}^d )} = \int_{F_2} dx \int_{\mathbf{R}^d} |g_h(x-y)|^2 |u(y)|^2 dy,
\label{glemma}  
\end{align}
where $g_h(x)$ takes the following form:
\begin{equation*}
g_h(x) = \left( \frac{1}{\pi h} \right)^{\frac{d}{4}} e^{-\frac{x^2}{2h}}=C_{d, h}e^{-\frac{x^2}{2h}} .
\end{equation*}
Then, by the assumption (\ref{prop5.2}), the right--hand side of (\ref{glemma}) is
\begin{align}
C_{d, h}^{2}  \int_{F_2} dx \int_{\mathbf{R}^d} e^{-\frac{(x-y)^2}{2h}}& e^{-\frac{(x-y)^2}{2h}} |u(y)|^2 dy \notag \\ 
& \le C_{d, h}^{2} e^{-\frac{\sigma_1^2}{2h}} \int_{F_2} dx \int_{\mathbf{R}^d} e^{-\frac{(x-y)^2}{2h}}|u(y)|^2 dy. \notag
\end{align}
Therefore 
\begin{equation}
||T u (\ \cdot \ , \ \cdot \ ; h)||^2_{L^2(F_2 \times \mathbf{R}^d )} \le C_{d, h}^{\frac{1}{2}} e^{-\frac{\sigma_1^2}{2h}} \int_{F_2} dx \int_{\mathbf{R}^d} e^{-\frac{(x-y)^2}{2h}}|u(y)|^2 dy.
\label{ippomae5.2}
\end{equation}
Since we have 
\begin{align*}
\int_{F_2} e^{-\frac{(x-y)^2}{2h}} dx \le \int_{\mathbf{R}^d} e^{-\frac{(x-y)^2}{2h}} dx = (2h)^{\frac{d}{2}} \int_{\mathbf{R}^d} e^{- \theta^2} d\theta = (2\pi h)^{\frac{d}{2}},
\end{align*}
(\ref{ippomae5.2}) implies the desired result 
\begin{align*}
||T u (\ \cdot \ , \ \cdot \ ; h)||^2_{L^2(F_2 \times \mathbf{R}^d )} \le (2 \pi h)^{\frac{d}{2}}C_{d, h}^2e^{- \frac{\sigma_1^2}{2h}} ||u||_{L^2(\mathbf{R}^d)}.
\end{align*}
\end{proof}
\end{prop}

\begin{prop}\label{prop-F-supp}
Let $F_1, F_2$ be a closed sets in $ \mathbf{R}^d $ satisfying 
\begin{equation}
{\rm dist}(F_1, F_2) = \sigma_{2} > 0.
\label{prop5.3supp}
\end{equation}
Then for $u \in L^2(\mathbf{R}^d)$ with $ {\rm supp} F_h u \subset F_1$, 
there exists $\delta > 0$ such that 
\begin{equation}
||T u (\ \cdot \ , \ \cdot \ ; h)||^2_{L^2(\mathbf{R}^d \times F_2)} \le 2^{\frac{d}{2}}e^{- \frac{\delta}{h}} ||u||_{L^2(\mathbf{R}^d)}.
\label{prop5.3final} 
\end{equation}
\begin{proof}
By Remark 3.4.4 in Martinez\cite{Mar}, we have 
\begin{align}
T u ( x , \xi ; h) = e^{i x \xi /h} T F_h u(\xi, -x).
\end{align}
Then by using Proposition \ref{prop-f-supp}, we are able to obtain easily the desired result.

\end{proof}
\end{prop}

\begin{thm}
Let $p=p(x,\xi;h)\in S_{2n}^{hol}(1,\Sigma_a)$, and set $P=Op_h^t(p)$ for $t\in [0,1]$ fixed.
Let $\psi=\psi(\xi)\in S_{n}(1)$ (independent of $h>0$) be real-valued, and assume
 $\sup_{\xi\in\mathbf{R}^n}|\nabla\psi(\xi)|<a$.
Let $f=f(x,\xi;h)\in S_{2n}(1)$. 
Then there exist $C>0$ and $h_0>0$ such that for every $u\in L^2(\mathbf{R}^n)$ and $h\in (0,h_0]$,
$$
\Big|\|fe^{\psi/h}TPu\|^2-\|f(x,\xi;h)p(x-i\nabla_{\xi}\psi(\xi), \xi-\nabla_{\xi}\psi(\xi);h)e^{\psi/h}Tu\|^2\Big|
\le Ch\|e^{\psi/h}Tu\|^2.
$$
\end{thm}

\begin{rem}
The estimate in the theorem above is uniform in $\varepsilon\in (0,1]$ 
if we replace $\psi$ with $\varepsilon \psi$.
\end{rem}

\section{Proof of main theorem}
In this section, using several propositions and lemmas in section 4 and section 6, we derive the estimation (\ref{outline}) and prove our main theorem.
\begin{proof}
We first write the integral equation (\ref{Ieq0}) as the sum of two parts as follows:
\begin{align}
\displaystyle&\int_{0}^{\tau^*}U_a(\tau^* - s) [w(s, \cdot) f(\cdot)](y) ds \notag \\
&=\int_{0}^{\tau^*_0} U_a(\tau^* - s) [w(s, \cdot) f(\cdot)](y) ds 
+ \int_{\tau^*_0}^{\tau^*} U_a(\tau^* - s) [w(s, \cdot) f(\cdot)](y) ds. \notag \\
& \equiv I_1(y) + I_2(y),
\label{mainI}
\end{align}
where $\tau_0^*$ is a positive constant such that $0<\tau_0^*<\tau^*.$
In the remaining part of this proof, to derive exponentially small $Tf$, we shall consider the $L^2$ estimate of (\ref{mainI}) with $H_a$, and we assume that $L_0 = L+1$.

First, by Lemma \ref{lem-I1}, we are able to get the estimation for $TH_aI_1(x, \xi;h)$ as (\ref{outline}).

Next, to consider the $L^2$ estimate of $TH_aI_2(x, \xi;h)$, we regard $H_a $ as a pseudodifferential operator acting on $f$(see \cite{Mar}), that is, 
\begin{align*}
H_aI_2(y) &= \int_{\tau^*_0}^{\tau^*} H_a U_a(\tau^* - s) [w(s, \cdot) f(\cdot)](y) ds \\
&={\rm Op}_h^1(\it{p})f(y),
\end{align*}
where the symbol of the above pseudodifferential operator takes the following form:
\begin{align}
p(y, \xi) &= \displaystyle(\xi + ia)^2 \int_{\tau^*_0}^{\tau^*} e^{-(\tau^*-s)(\xi+ia)^2} w(s, y) ds  \notag \\ 
&= \displaystyle \int_{\tau_0^*}^{\tau^*} \frac{\partial}{\partial s} (e^{-(\tau^*-s)(\xi+ia)^2}) w(s, y) ds \notag \\
&= w(\tau^*, y) - e^{-(\tau^*-\tau_0^*)(\xi+ia)^2}w(\tau_0^*, y) 
- \int_{\tau_0^*}^{\tau^*} e^{-(\tau^*-s)(\xi+ia)^2} \frac{\partial w}{\partial s} (s, y) ds.
\label{p-eq}
\end{align}
Here, let $\chi_1(\xi) \in C^\infty_0(\textbf{R})$ be such that $\chi_1=0$ if $|\xi|<\frac{1}{4}, \chi_1 = 1$ if $|\xi| > \frac{1}{2}$ and we set 
\begin{equation}
p_j(x, \xi; h) = p\left(x, \xi/h \right) \chi_j(\xi) \hspace{1ex} (j=1, 2),
\label{pi}
\end{equation}
where $\chi_2(\xi) = 1 - \chi_1(\xi).$
Moreover, let the real-valued function $\psi \in C_0^{\infty}(\textbf{R})$ be such that $\psi = 0$ if $|\xi| < 1, \psi = 1$ if $|\xi| > 2$ and there exist $\varepsilon_0 > 0$ such that $\varepsilon_0 ||\nabla \psi||_{L^{\infty}} < \sigma_0$, where $\sigma_0$ is a constant in Lemma \ref{lem-w}. 

Now we apply Corollary 3.5.5 (in \cite{Mar}) with $T=T^{\varepsilon}, f=1$ and $\psi = \varepsilon \psi$, where $\varepsilon > 0$ will be taken small enough later and we set
\begin{equation}
T^{\varepsilon} u  = e^{\varepsilon \psi(\xi) / h} T u. 
\label{T-e}
\end{equation}
Then, we obtain  
\begin{align}
||T^{\varepsilon}&{\rm Op}_h^1(p_1)f||_{L^2}^2 \notag \\
&\ge ||p_1(y-i\varepsilon \partial_{\xi} \psi(\xi), \xi - \varepsilon \partial_{\xi} \psi(\xi) ; h) T^{\varepsilon} f||^2_{L^2}  -C h||T^{\varepsilon} f||_{L^2}^2. 
\label{mar11}
\end{align}
Using Taylor's formula and Lemma \ref{lem-p1}, we can estimate the right--hand side of (\ref{mar11}) as follows:
\begin{align*}
&\ge ||p_1(y, \xi ; h) T^{\varepsilon} f||^2_{L^2} - C_1(\varepsilon + h)||T^{\varepsilon} f||^2_{L^2} \notag \\
\mbox{} \\
&\ge ||w(\tau^*, y)\chi_1 T^{\varepsilon} f||^2_{L^2} - C_1(\varepsilon + h)||T^{\varepsilon} f||^2_{L^2}. \notag \\
\end{align*}
By $\chi_1 = 1 - \chi_2$,
\begin{align*}
||w(\tau^*, y)(1-\chi_2) T^{\varepsilon} f||^2_{L^2}& - C_1(\varepsilon + h)||T^{\varepsilon} f||^2_{L^2} \notag \\
\ge \frac{1}{2}||w(\tau^*, y) &T^{\varepsilon} f||_{L^2}^2 - 2||w(\tau^*, y) \chi_2 T^{\varepsilon} f||_{L^2}^2 
-C_1(\varepsilon + h)||T^{\varepsilon} f||^2_{L^2} \notag\\
\end{align*}
and (\ref{T-e}) and Lemma \ref{lem-w} we have
\begin{align*}
\ge \left\{ \frac{C_0^2}{2}-C_1(\varepsilon + h) \right\}  ||T^{\varepsilon} &f||^2_{L^2([-L_0, \infty) \times \textbf{R})} \notag\\
-C^2 ||Tf||^2_{L^2}& -C_1(\varepsilon + h) ||T^{\varepsilon} f||^2_{L^2((-\infty, -L_0) \times \textbf{R})}
\end{align*}
where we used $\psi(\xi) = 0$ if $|\xi| \le 1.$ 

On the other hand, since we can rewrite the left--hand side of (\ref{mar11}) as
\begin{align*}
||T^{\varepsilon}{\rm Op}_h^1(p_1) f||_{L^2}^2 &= ||T^{\varepsilon}[-H_aI_1(y) - {\rm Op}_h^1(p_2) f]||_{L^2}^2 \notag \\ & \le 2||T^\varepsilon H_a I_1(y)||_{L^2}^2 + 2||T^{\varepsilon}{\rm Op}_h^1(p_2) f||_{L^2}^2,
\end{align*}
by using the following
\begin{align*}
||T^{\varepsilon}H_a I_1||^2_{L^2(\mathbf{R} \times \{|\xi| \le 1 \})}
= ||T H_a I_1||^2_{L^2(\mathbf{R} \times \{|\xi| \le 1 \})} \le C_2 ||Tf||^2_{L^2}
 \notag
\end{align*}
and
\begin{align*}
||T^{\varepsilon}&{\rm Op}_h^1(p_2)f||^2_{L^2(\mathbf{R} \times \{|\xi| \le 1 \})}
= &||T {\rm Op}_h^1(p_2)f||^2_{L^2(\mathbf{R} \times \{|\xi| \le 1 \})}
\le C_3 ||Tf||^2_{L^2},
\notag
\end{align*}
\noindent we have the following estimates 
\begin{align}
||T^{\varepsilon}&{\rm Op}_h^1(p_1) f||^2_{L^2} \notag \\
&\le ||T^{\varepsilon}H_a I_1||^2_{L^2(\mathbf{R} \times \{|\xi| \le 1 \})}
 + ||T^{\varepsilon}H_a I_1||^2_{L^2(\mathbf{R} \times \{|\xi| > 1 \})}
 \notag \\
&\hspace{7ex} + ||T^\varepsilon {\rm Op}_h^1(p_2) f]||^2_{L^2(\mathbf{R} \times \{|\xi| \le 1 \})}
 + ||T^\varepsilon {\rm Op}_h^1(p_2) f]||^2_{L^2(\mathbf{R} \times \{|\xi| > 1 \})}
  \notag \\ 
& \le C_4||T f||_{L^2}^2 \notag \\
&\hspace{7ex} + ||T^{\varepsilon}H_a I_1||^2_{L^2(\mathbf{R} \times \{|\xi| > 1 \})}
+||T^\varepsilon {\rm Op}_h^1(p_2) f]||^2_{L^2(\mathbf{R} \times \{|\xi| > 1 \})}
.
\label{lhside}
\end{align}
Then, these estimates give 
\begin{align*}
\left\{\frac{C_0^2}{2} -C_1  (\varepsilon + h) \right\} & ||T^{\varepsilon} f||^2_{L^2([-L_0, \infty) \times \textbf{R})} \notag\\
 -C^2 ||Tf&||^2_{L^2}
 -C_1(\varepsilon + h) ||T^{\varepsilon} f||^2_{L^2((-\infty, -L_0) \times \textbf{R})} \notag \\
\le C_4||T f ||_{L^2}^2 &+ ||T^{\varepsilon}H_a I_1||^2_{L^2(\mathbf{R} \times \{|\xi| > 1 \})}+||T^\varepsilon {\rm Op}_h^1(p_2) f]||^2_{L^2(\mathbf{R} \times \{|\xi| > 1 \})}.
\end{align*}
Since we can get the following by applying to Proposition \ref{prop-f-supp}, Lemma \ref{lem-I1} and Lemma \ref{lem-p2}
\begin{align*}
||T^{\varepsilon} f||^2_{L^2((-\infty, -L_0) \times \textbf{R})} 
= ||e^{\varepsilon \frac{\psi(\xi)}{h}}Tf||^2_{L^2((-\infty, -L_0) \times \textbf{R})} =  O(e^{\frac{2\varepsilon-\delta_1}{h}})||Tf||^2_{L^2},
\end{align*}
\begin{align*}
||T^{\varepsilon}H_a I_1||^2_{L^2(\mathbf{R} \times |\xi| > 1)}
= ||e^{\varepsilon \frac{\psi(\xi)}{h}}H_a I_1||^2_{L^2(\mathbf{R} \times \{|\xi| > 1 \})}=  O(e^{\frac{2\varepsilon - \delta_2}{h} })||Tf||^2_{L^2}
\end{align*}
and
\begin{align*}
||T^\varepsilon {\rm Op}_h^1(p_2) f]||^2_{\mathbf{R} \times L^2(|\xi| > 1)}
=||e^{\varepsilon \frac{\psi(\xi)}{h}}{\rm Op}_h^1(p_2) f]||^2_{L^2(\mathbf{R} \times \{|\xi| > 1 \})}
=  O(e^{\frac{2\varepsilon - \delta_3}{h}})||Tf||^2_{L^2},
\end{align*}
where $\delta_i \ (i=1,2,3)>0$ are some constant.
Then, for $\delta_i>0$ if $\varepsilon$ is chosen small enough, 
we have 
\begin{align*}
||e^{\varepsilon \frac{\psi(\xi)}{h}}T f ||^2_{L^2 ([-L_0, \infty) \times \textbf{R})} = O(1) ||T f ||_{L^2}^2.
\label{ketsuronmae}
\end{align*}
Since $\psi(\xi) = 0$ if  $|\xi| \ge 2, $ we obtain
\begin{equation}
||Tf||^2_{L^2 ([-L_0, \infty) \times \{ |\xi| \ge 2 \})} = O(e^{-\frac{\delta}{h}}).
\label{ketsuron}
\end{equation}
In particular, we deduce from (\ref{ketsuron}) that 
\begin{equation*}
[-L_0, \infty) \times \{ |\xi| \ge 2 \} \cap {\rm WF}_a(f) = \emptyset.
\end{equation*}
Hence, we obtain that $f$ is real analytic in $(-L_0, \infty)$.

On the other hand, since $f=0$ in $(-L_0, -L)$ by the assumption, we conclude that $f$ is identically zero on $\textbf{R}$.

The proof is complete.

\end{proof}

\section{Lemmas}
In this section, we prove some of the auxiliary mathematical results which plays an 
essential role in the proof of main theorem.
First, we define the following functional spaces $L^2_\lambda, H^s_\lambda$:
\begin{eqnarray*}
L^2_\lambda(\mathbf{R}) = \{ u \in L^2_{loc}(\mathbf{R}); \hspace{2ex} e^{\lambda <x>} u \in L^2(\mathbf{R}) \} \hspace{2ex} (\lambda \in \mathbf{R}) \\
H^s_\lambda(\mathbf{R}) = \{ u \in D'(\mathbf{R}); \hspace{2.5ex} e^{\lambda <x>} u \in H^s(\mathbf{R}) \} \hspace{2ex} (\lambda \in \mathbf{R}).
\end{eqnarray*}
Then we can prove the following result concerning the direct problem (\ref{mueq3}).
\begin{lem}\label{lem-direct}
Assume that $f \in L^2( [0, T] ; H_\lambda^{s-1}(\mathbf{R}) ).$ 
Then there is a solution $u \in C([0, T] ; H_\lambda^s(\mathbf{R})) \cap L^2([0, T] ; H_\lambda^{s+1}(\mathbf{R}))$;  moreover, the solution is unique in $\displaystyle \cup_{s', \lambda'} L^2([0, T] ; H_{\lambda'}^{s'}(\mathbf{R})).$   
\begin{proof}
This result is found in the book by Friedman \cite{Fri}.
\end{proof}
\end{lem}

The properties of $w$ are as follows:

\begin{lem}\label{lem-w}

\mbox{}

(i) For any $\tau > 0$ and $ y \in \mathbf{R}$,  $|w(\tau, y)| \le e^{a^2 \tau}. $

(ii) $w(\tau, y) \in C^{\infty}((0, \infty) \times \mathbf{C}_z)$. Moreover, for $\tau>0$, $w(\tau, \ \cdot \ )$ can be extended as a holomorphic function of $z$ on $\mathbf{C}$.

(iii) For $\sigma_0 > 0, k = 0, 1, 2, \cdots, \ \alpha = 0, 1, 2, \cdots$ and $\tau_0$ such that $0<\tau_0<\tau^*$, there exists $C>0$ such that 
\begin{equation*}
|\partial_{\tau}^{\alpha}\partial_{z}^{\alpha}w(\tau, z)| \le C \hspace{5ex}  \mbox{for any} \tau \in [\tau_0, \tau^*] \  \mbox{and}  z \in \mathbf{C} \ \mbox{such that} \ |Im z| \le \sigma_0,
\end{equation*}
where $C$ depends on $\sigma_0, k, \alpha$ and $\tau_0.$

(iv) For $\tau_0$ such that $0 < \tau_0 < \tau^* $ and $L_0 \ge 0$, there exists $C_0 > 0$ such  that
\begin{equation} 
w(\tau, y) \ge C_0 \hspace{5ex} \mbox{for any} \  \tau \in [\tau_0, \tau^*] \ \mbox{and} \  y \ge -L_0.
\end{equation}

\begin{proof}
Assertions (i)$\sim$(iv) easily follow from the form of $w$ and Cauchy$'$s integral formula.
\end{proof}
\end{lem}
\begin{lem}\label{lem-I1}
For all $C>0$, there exists $\delta > 0$ such that
\begin{equation}
\lVert TH_aI_1(\ \cdot, \ \cdot \ ; h) \rVert^2_{L^2(\textbf{\mbox{R}} \times \{|\xi| \ge C \})} = O (e^{-\frac{\delta}{h}}).
\end{equation}
\begin{proof}
Since we transform $TH_a I_1(y)$ into the following form,  
\begin{align*}
TH_a I_1(y) &= TH_a U_a(\tau^*-\tau^*_0) \displaystyle \int_0^{\tau_0^*} U_a(\tau_0^* - s) [ w(s, \cdot) f(\cdot)] (y) ds \notag \\
\end{align*}
To prove the assertion we have only to prove that if $f \in L^2(\textbf{R})$, then
\begin{equation}
\int_0^{\tau_0^*} U_a(\tau_0^* - s) [ w(s, \cdot) f(\cdot)] (y) ds 
\label{intw}
\end{equation}
belongs to $L^2(\textbf{R}).$

Now we set 
\begin{equation*}
f_1(y) := \int_0^{\tau_0^*} U_a(\tau_0^* - s) [ w(s, \cdot) f(\cdot)] (y) ds. 
\end{equation*}
Then we have
\begin{equation*}
||f_1||_{L^2} \le \int_0^{\tau_0^*}|| U_a(\tau_0^* - s) [ w(s, \cdot) f(\cdot)](\cdot)||_{L^2} ds 
\le \int_0^{\tau_0^*}|| w(s, \cdot) f(\cdot)||_{L^2} ds
\end{equation*}
By Lemmma \ref{lem-w} (i), 
$w(s, x) f(x)$ belongs to $L^{\infty}([0, \tau^*_0] \ ; \ L^2(\textbf{R})).$
\noindent Hence we can show that the integral (\ref{intw}) belongs to $L^2(\textbf{R})$.

Next, since $F_h H_a I_1(\xi)$ is the following form
\begin{equation}
F_h  H_a I_1(\xi) = \left(\frac{\xi}{h} + ia\right)^2 e^{-(\tau^* - \tau_0^* ) \left(\frac{\xi}{h} + ia\right)^2} F_h f_1(\xi),
\end{equation}
we have
\begin{align}
\displaystyle \lVert TH_aI_1(\ \cdot, \ \cdot \ & ; h)\rVert^2_{L^2(\textbf{\mbox{R}} \times \{ \lvert \xi \rvert \ge C\} )} \notag \\
&=\displaystyle \int_{\textbf{R}} \int_{\lvert \xi \rvert \ge C} \lvert g_h(\xi - \eta) \rvert^2 \lvert F_h H_a I_1(\eta) \rvert^2 d\xi d\eta \notag \\
&\le \displaystyle \frac{1}{(\pi h)^{\frac{1}{2}}} \int_{\textbf{R}} \int_{\lvert \xi \rvert \ge C} \left(\frac{\xi}{h} + ia\right)^4  e^{-\delta_1 \frac{\lvert \xi - \eta \rvert^2}{h}}e^{-\delta_1 \frac{\lvert \eta \rvert^2}{h}} \lvert F_h f_1\left(\eta \right) \rvert^2 d\xi d\eta \notag \\
&\le \displaystyle C_1 \int_{\textbf{R}} \int_{\lvert \xi \rvert \ge C}\left(\frac{\xi}{h} + ia\right)^4 e^{-\delta_2 \frac{\lvert \eta \rvert^2 + \lvert \xi \rvert^2}{h}} \lvert F_h f_1\left(\eta\right) \rvert^2 d\xi d\eta \notag \\
&\le \displaystyle C_1 e^{-\delta_2\frac{C }{2 h}}  \int_{\textbf{R}} \int_{\lvert \xi \rvert \ge C} \left(\frac{\xi}{h} + ia\right)^4 e^{-\delta_2 \frac{\lvert \xi \rvert^2}{2h}} e^{-\delta_2 \frac{\lvert \eta \rvert^2}{h}} \lvert F_h f_1\left(\eta\right) \rvert^2 d\xi d\eta \notag.
\end{align}
Therefore we obtain 
\begin{align*} 
\lVert TH_aI_1(\ \cdot, \ \cdot \ ; h)\rVert^2_{L^2(\textbf{\mbox{R}} \times \{ |\xi| > C\}) } \le C_2 e^{-\frac{\delta}{h}} ||f_1||_{L^2(\mathbf{R})}.
\end{align*}

The proof is complete.
\end{proof}
\end{lem}

\begin{lem}\label{lem-p}
\mbox{}

\noindent (i) For $\alpha, \beta \in \textbf{N}$, there exists $C > 0$ for any $z$ such that $|Im z| < \sigma_0, \xi \in \mathbb{R}$, 
\begin{equation}
|\partial_z^{\alpha} \partial_\xi^\beta p(z, \xi) | = O(\langle \xi \rangle^{-|\beta|}),
\end{equation}
where $\sigma_0$ is a constant in Lemma \ref{lem-w}.
 
\noindent (ii) There exists $C > 0$ such that for any $x, \xi \in \textbf{R}$,
\begin{equation}
|p(x, \xi) - w(\tau^*, x)| = O(\langle \xi \rangle^{-2}).
\end{equation} 
\begin{proof}
Assertion (i) is obtained immediately by (\ref{p-eq}) and Lemma \ref{lem-w}.

For (ii), we rewrite (\ref{p-eq}) as     
\begin{align}
p(y, \xi) - w(\tau^*, y) = -e^{-(\tau^*-\tau_0^*)(\xi+ia)^2}w(\tau_0^*, y) 
- \int_{\tau_0^*}^{\tau^*} e^{-(\tau^*-s)(\xi+ia)^2} \frac{\partial w}{\partial s} (s, y) ds.
\label{p-eq2}
\end{align}
Since the first term on the right--hand side of (\ref{p-eq2}) is $O(e^{-|\xi|^2})$ and the second integral on the right--hand side of (\ref{p-eq2}) is $O( \langle \xi \rangle^{-2})$, we can obtain assertion (ii). 

\end{proof}
\end{lem}

Here we prove the following lemma for $p_1(x, \xi ; h).$

\begin{lem}\label{lem-p1}
\mbox{}

(i) $p_1(x, \xi; h) \in S_2^{hol}(1, \Sigma_{\sigma_0})$.

(ii) $|p_1(x, \xi; h) - w(\tau^*, x) \chi_1(\xi)| = O(h^2)$.

\begin{proof}
Assertion (i) is obtained immediately by (\ref{p-eq}) and the definition of $p_1(x, \xi; h)$ in (\ref{pi}).

For (ii), by  Lemma \ref{lem-p} (ii) we have
\begin{equation}
\lvert p(x, \xi/h) - w(\tau^*, x) \rvert \le C \langle \xi/h \rangle^{-2}.
\end{equation} 
Since $p_1(x, \xi; h) = \chi_1(\xi)$, we have
\begin{align*}
\lvert p_1(x, \xi/h) -  \chi_1&(\xi)w(\tau^*, x) \rvert  
\le C \chi_1(\xi) \langle \xi/h \rangle^{-2} 
\le C_0 h^2,
\end{align*}
where $C_0$ is a constant independent of $\xi$.

Thus the proof is complete.
\end{proof}
\end{lem}

\begin{lem}\label{lem-p2}
Let $p_2 = p_2(x, \xi; h)$, then there exists $\delta > 0$ such that 
\begin{equation}
\lVert T[{\rm OP}_{h}^1(p_2)f]( \ \cdot \ , \ \cdot \ ; h)\rVert^2_{L^2(\mathbf{R} \times \{ |\xi| \ge \frac{3}{4} \})} = O(e^{- \frac{\delta}{h}} ) \lVert Tf\rVert_{L^2}^2.
\end{equation}
\begin{proof}
By definition
\begin{align}
{\rm OP}_{h}^1(p_2)&f(x) \notag \\
&= (2 \pi h)^{-1} \displaystyle \int_{\mathbf{R}} \int_{\mathbf{R}} e^{i(x-y) \xi /h} p_2(y, \xi ; h) f(y) dy d\xi  \notag \\
\end{align}
which implies 
\begin{equation*}
F_h[{\rm OP}^1_h(p_2)f(\cdot)](\xi) = (2 \pi h)^{-\frac{1}{2}} \displaystyle \int_{\mathbf{R}}  e^{-i y\xi/h} p_2(y, \xi ; h) f(y) dy.
\end{equation*}
Hence by (\ref{pi}), we have
\begin{equation*}
{\rm supp} F_h[{\rm OP}^1_h(p_2)f(\cdot)](\xi) \subset \left\{ |\xi| \le 1/4 \right\}.  
\end{equation*}
Therefore from Proposition \ref{prop-F-supp} we can obtain the conclusion.

The proof is complete.

\end{proof}
\end{lem}


\begin{thebibliography}{33}


\bibitem{B-S} Black F and Scholes M. \textit{The pricing of options and corporate liabilities.} Journal of Political Economy, 1973, \textbf{81}, 637-659. 

\bibitem{Boness} Boness A.  \textit{Elements of a theory of stock-option value.} J. Political Econ, 1964, \textbf{72}(2), 163-175. 

\bibitem{B-I} Bouchouev I and Isakov V. \textit{Uniqueness, stability and numerical methods for the inverse problem that arises in financial markets.} Inverse Problems, 1999,  \textbf{15}, R95-R116.

\bibitem{B-I-V}

Bouchouev I , Isakov V and Valdivia N. \textit{Recovery of volatility coefficient by linearrization.} Quantitative Finance, 2002, Vol2, 257-263.

\bibitem{dupire}

Dupire B. \textit{Pricing with a smile.} Risk 7, 1994, 18-20.

\bibitem{Fri} Friedman A. \textit{Partial Differential Equations of Parabolic Type.} 1964, Prentice-Hall, Englewood Cliffs, N.J.

\bibitem{kubo} Korolev M, Kubo H and Yagola G. \textit{Parameter identification problem for a parabolic equation-application to the Black-Scholes option pricing model}2012, J. Inverse Ill-posed probl. \textbf{20} No.3 327-337 

\bibitem{Mar} Martinez A. \textit{An introduction to Semiclassical and Microlocal Analysis.} 2002, Universitext Springer.

\bibitem{merton}
Merton R. \textit{On the Pricing of Corporate Debt: The Risk Structure of Interest Rates.} Journal of Finance, 29, 1974, 449-470.

\bibitem{myy} Mitsuhiro M, Ota Y and Yadohisa H \textit{New mathematical approach for an inverse problem in financial markets. } Proceedings of COMPSTAT2012, 2012, 585-594.

  




\end{thebibliography}
\end{document}